\documentclass[a4paper,draft]{amsart}
\usepackage[a4paper,margin=25mm]{geometry}
\usepackage{amsmath}
\usepackage{amssymb}
\usepackage{url}

\parskip\smallskipamount

\begin{document}

 \author[Manuel Kauers]{Manuel Kauers\,$^\ast$}
 \address{Manuel Kauers, Institute for Algebra, J. Kepler University Linz, Austria}
 \email{manuel.kauers@jku.at}

 \author[Martina Seidl]{Martina Seidl\,$^{\ast\ast}$}
 \address{Martina Seidl, Institute for Formal Models and Verification, J. Kepler University Linz, Austria}
 \email{martina.seidl@jku.at}

 \author[Doron Zeilberger]{Doron Zeilberger}
 \address{Doron Zeilberger, Department of Mathematics, Rutgers University, USA}
 \email{DoronZeil@gmail.com}
 
 \thanks{$^\ast$ Supported by the Austrian FWF grants F5004 and P31571-N32.}

 \thanks{$^{\ast\ast}$ Supported by the Austrian FWF grant S11408-N23.}
 
 \title{On the maximal minimal cube lengths in distinct~DNF~tautologies}

 \begin{abstract}
   Inspired by a recent article by Anthony Zaleski and Doron Zeilberger, we investigate
  the question of determining the largest $k$ 
  for which there exists  boolean formulas in disjunctive normal form (DNF) with $n$ variables,
  none of whose conjunctions are `parallel',
  and   such that  all of them have at least $k$ literals. 
 Using a SAT solver, we answer some of the questions they left open.
   We also determine the corresponding numbers for DNFs obeying certain symmetries.
 \end{abstract}

 \maketitle

 \section{Problem Statement}

 We consider boolean formulas with $n$ variables $x_1,\dots,x_n$.
 A \emph{literal} is a variable or a negated variable, e.g., $x_3$ or $\bar x_7$.
 A \emph{cube} is a conjunction of literals, e.g., $x_3\land\bar x_7$.
 The length of a cube is the number of distinct literals appearing in it.
 A formula in \emph{disjunctive normal form} (DNF) is a disjunction of cubes, e.g.,
 $(x_3\land\bar x_7)\lor(x_5\land\bar x_6\land x_7)$.
 Such a DNF is called a \emph{tautology} if it evaluates to true for all
 assignments of the variables.
 For example $x_3\lor x_5\lor(\bar x_3\land\bar x_5)$ is a tautology.
 It consists of two cubes of length~1 and one cube of length~2.

 Inspired by work of Erd\"os~\cite{erdos},
 Zaleski and Zeilberger~\cite{zz} have recently considered DNFs in which
 all cubes have distinct supports. 
The support of a cube is the set of  variables occurring in it. 
 For example the support of the cube $x_3$ is the singleton set $\{3\}$,
the support of the cube $\bar x_5$ is the singleton set $\{5\}$,
while the support of the cube $\bar x_3\land x_5$ is the set $\{3,5\}$.
This implies that the DNF $x_3\lor {\bar x_5} \lor (\bar x_3\land x_5)$  has distinct supports. 
On the other hand the Hamlet question $x_1 \lor \bar x_1$ does not have distinct supports.
They call these formulas \emph{distinct DNFs.}
 Inspired by a study of covering systems, Zaleski and Zeilberger want to know,
 for any given~$n$, what is the largest $k$ such that there is a distinct DNF
 tautology with $n$ variables only consisting of cubes of length at least~$k$.

 Using a greedy algorithm, they searched for distinct DNF tautologies with
 a prescribed number of variables and a prescribed minimal cube length. The
 largest minimal cube length for which they found formulas are as follows:
 \begin{center}
   \begin{tabular}{c|cccccccccccccccccc}
     $n$ & $1$ & $2$ & $3$ & $4$ & $5$ & $6$ & $7$ & $8$ & $9$ & $10$ & $11$ & $12$ & $13$ & $14$  \\\hline
     $k$ & $0$ & $1$ & $1$ & $2$ & $3$ & $4$ & $4$ & $5$ & $6$ & $6$\rlap{$_?$} & $7$ & $8$ & $9$ & $9$\rlap{$_?$}
   \end{tabular}
 \end{center}
 These are only lower bounds for the optimal values of~$k$. However, by a density
 argument it can be shown that the optimal $k$ must satisfy the inequality $\sum_{i=k}^n\binom ni2^{-i}\geq1$,
 which gives rise to upper bounds.
 The numbers given in the table above turn out to match the upper bounds except for $n=10$ and $n=14$
 (indicated by  question marks), where they are off by one.

 As a variant of the problem, Zaleski and Zeilberger also wanted to know, for
 any given $n$, what is the largest $k$  such that there is a distinct DNF
 tautology with $n$ variables only consisting of cubes of length exactly~$k$.
 In this case, the density argument implies that such a $k$ must satisfy $\binom nk2^{-k}\geq1$, which
 again gives an upper bound. With their greedy approach, they determined the
 following lower bounds. Again, mismatches with the upper bound are indicated
 by a question mark.
 \begin{center}
   \begin{tabular}{c|cccccccccccccccccc}
     $n$ & $1$ & $2$ & $3$ & $4$ & $5$ & $6$ & $7$ & $8$ & $9$ & $10$ & $11$ & $12$ & $13$ & $14$  \\\hline
     $k$ & $0$ & $0$ & $0$\rlap{$_?$} & $2$ & $2$\rlap{$_?$} & $3$ & $4$ & $5$ & $5$\rlap{$_?$} & $6$ & $7$ & $8$ & $8$\rlap{$_?$} & $9$
   \end{tabular}
 \end{center}
 It is clear that there is no solution for $n=3$ and $k=1$, so in this case the upper bound is too pessimistic
 and $k=0$ is the right value.

 For $n=5$ and $n=9$ the computations reported in the present paper imply that
 the values $2$ and $5$ are also correct.
 We were not able to confirm the entry for $n=10$ in the first table with about one year of computation time. 
 We did not attempt
 to confirm the entries for $n=14$ in the first or $n=13$ in the second table.

 We add two refinements to the problem. First, we introduce an additional parameter~$u$
 which bounds the lengths of the cubes from above. For any particular choice $n,u$,
 we want to know the largest $k$ such that there is a distinct DNF tautology with
 $n$ variables only consisting of cubes of length at least $k$ and at most~$u$.
 The special case $u=n$ corresponds to the first variant of Zaleski and Zeilberger
 and the special case $u=k$ corresponds to the second variant. We think that the
 intermediate cases are also of interest.

 Our second refinement concerns symmetries. Letting permutations act on the indices
 of the variables, we say that a DNF is \emph{invariant} under a certain subgroup $G$
 of~$S_n$ if every $g\in G$ maps the DNF to itself. For example, the
 DNF $(x_1\land\bar x_2\land x_3)\lor
 (x_2\land\bar x_3\land x_4)\lor
 (x_3\land\bar x_4\land x_1)\lor
 (x_4\land\bar x_1\land x_3)$ is invariant under the cyclic group~$C_4$.
 For the groups $C_n$,~$D_n$, $A_n$, and~$S_n$, and for various choices of $n$ and~$u$,
 we have determined the largest $k$ such that there is a distinct DNF tautology with $n$
 variables consisting of cubes of lengths at least $k$ and at most~$u$ which are invariant
 under the given group. 
 
 \section{SAT Encoding}

 Our results were obtained with the help of a SAT solver~\cite{heule,knuth}, using a rather straightforward
 encoding of the problem. For each cube, we introduced one boolean variable that indicates
 whether or not this cube is going to be a part of the DNF we are looking for. Note that
 this creates $\sum_{i=k}^u\binom ni2^i$ variables, a quantity that grows quickly when $n$
 or $u-k$ increase. For example, in the case $n=u=10$ and $k=7$, where we were unable to
 complete the computation, we were dealing with 33024 variables.

 In order to enforce that the DNF is a tautology, we specify for every
 assignment a clause saying that at least one of the cubes that becomes true
 under this assignment must be selected. In order to enforce that the DNF be
 distinct, we have to specify clauses which encode the requirement that for
 every possible support, at most one of the cubes having this support can be
 selected. There are many ways to encode a constraint of the form ``at most
 one'', and their pros and cons are discussed extensively in the
 literature~\cite{<=1,<=1b}. For our purpose, the so-called binary encoding
 seemed to work well.

 Finally, in order to enforce invariance under a certain group, we chose a set of generators
 and added for each cube $c$ and each generator $g$ a clause that says ``if $c$ is selected,
 then also $g(c)$''.


 The encoding as described so far is sufficient for proving existence or
 non-existence of a distinct DNF tautology for any prescribed $n,u,k$, and any
 prescribed group. In order to speed up the computations in practice, we may add
 some further constraints. One idea is to add clauses which forbid to select two
 cubes where one is strictly contained in the other. This is clearly a valid
 restriction, because when there is a solution that has two cubes that are
 contained in one another, we can discard the smaller one from it and obtain
 another solution.  However, it turns out that this particular idea floods the
 formula with too many additional clauses and slows down the computation rather
 than speeding it up.

 It is more efficient to break the symmetry of the problem, a standard
 technique in the context of SAT solving~\cite{sym}. Clearly, when there is a
 distinct DNF for certain $n,u,k$ and a certain group, then permuting all the
 variables $x_1,\dots,x_n$ in some way will yield another solution. Also
 replacing a certain variable $x_i$ by its negation $\bar x_i$ (and canceling
 double negation introduced by that) turns a solution into a new solution.
 Since we dropped the idea to forbid cubes that are contained in other cubes, we
 can restrict the search to a solution containing a cube of length~$k$, and
 because we are free to permute and negate variables, we may assume this
 cube to be $x_1\land x_2\land\cdots\land x_k$.

 Adding the variable for this cube to the formula allows for an appreciable
 amount of simplification (called unit propagation~\cite{heule} in SAT jargon). We are left
 with the freedom to permute the variables $x_1,\dots,x_k$ and the variables
 $x_{k+1},\dots,x_n$. By the first, it is fair to enforce an assumption that the
 variables are indexed in such a way that when a cube with support $x_1,\dots,x_{k-1},x_{k+1}$
 is selected, there is some $i$ such that $x_1,\dots,x_i$ appear negated in it
 and the remaining variables do not. This assumption may still leave some degrees
 of freedom, which can be used to make a similar restriction as to which cubes
 with support $x_1,\dots,x_{k-1},x_{k+2}$ may be selected.
 The freedom to permute the variables $x_{k+1},\dots,x_n$ is exploited by
 restricting the search to DNFs such for every $i=k+1,\dots,n-1$,
 the cube $x_1\land\cdots\land x_{k-1}\land x_{i+1}$ is only selected when 
 $x_1\land\cdots\land x_{k-1}\land x_i$ is also selected. 

 \section{Results}

 We have written a Python script that produces the SAT instances described in the
 previous section, and we have used Biere's award-winning SAT solver Treengeling~\cite{armin} to
 solve them. The results are summarized in the following tables, in which $n$ appears
 increased towards the right and $u$ grows downwards. Entries with $u>n$ are left blank
 because they are equivalent to $u=n$.

 By the density argument, the maximal $k$ for a particular choice of $n$ and $u$ must
 satisfy the inequality $\sum_{i=k}^u\binom ni2^{-i}\geq1$. In the following table,
 an entry is boxed when it does not match this bound. For the entries marked with
 a question mark, we have not been able to prove that the $k$ we found is really optimal,
 but the long and successless search is at least some indication that the bound is not
 reached in these cases. For $(n,u)\in\{(5,3),(9,6),(10,8)\}$, the SAT solver is able
 to show that distinct DNF tautologies with $k=3$, $k=6$, $k=7$, respectively, do not
 exist, although their existence would not be in conflict with the density bound.

 \fboxsep=1.5pt
 \def\mark#1{\smash{\hbox to0pt{\hss\fbox{#1}\hss}}}
 
    \begin{center}\small
      \begin{tabular}{c|cccccccccc}
       & 2 & 3 & 4 & 5 & 6 & 7 & 8 & 9 & 10 \\\hline
     2 & 1 & 1 & 2 & 2 & 2 & 2 & 2 & 2 & 2 \\
     3 &   & 1 & 2 & \mark{2} & 3 & 3 & 3 & 3 & 3 \\
     4 &   &   & 2 & 3 & 3 & 4 & 4 & 4 & 4 \\
     5 &   &   &   & 3 & 4 & 4 & 5 & 5 & 5 \\
     6 &   &   &   &   & 4 & 4 & 5 & \mark{5} & 6 \\
     7 &   &   &   &   &   & 4 & 5 & 6 & 6 \\
     8 &   &   &   &   &   &   & 5 & 6 & \mark6 \\
     9 &   &   &   &   &   &   &   & 6 & \mark6\rlap{\kern5pt$_?$} \\
     10 &   &   &   &   &   &   &   &    & \mark6\rlap{\kern5pt$_?$} 
  \end{tabular}
  \end{center}

  The next two tables contain our results about distinct DNF tautologies invariant
  under certain groups. We have investigated the cyclic group $C_n$, the
  dihedral group $D_n$, the alternating group $A_n$, and the full symmetric
  group~$S_n$. The table on the left lists the numbers for $C_n$ and~$D_n$,
  which turn out to be identical. Boxed entries highlight the differences to the
  previous table. The question marks refer to the search for~$C_n$, which for
  three entries did not terminate in a reasonable amount of time. Interestingly,
  it follows from the previous table that the entry for $(n,u)=(10,8)$ is 
  correct, but while the SAT solver was able to prove this in the (seemingly
  harder) case without invariant constraints, it did not succeed with the
  constraints for~$C_n$. The computations for all entries terminated in presence
  of the constraints for~$D_n$.
  
   \begin{center}\small
        \begin{tabular}{c|cccccccccc}
            & 2 & 3 & 4 & 5 & 6 & 7 & 8 & 9 & 10 \\\hline
          2 & 1 & 1 & 1 & 1 & 1 & 1 & 1 & 1 & 1 \\
          3 &   & 1 & 2 & 2 & 2 & 3 & 3 & 3 & 3 \\
          4 &   &   & 2 & \mark2 & 3 & \mark3 & 4 & 4 & 4 \\
          5 &   &   &   & \mark2 & \mark3 & 4 & \mark4 & 5 & 5 \\
          6 &   &   &   &   & \mark3 & 4 & 5 & 5 & 6 \\
          7 &   &   &   &   &   & 4 & 5 & 6 & 6 \\
          8 &   &   &   &   &   &   & 5 & 6 & 6\rlap{$_?$} \\
          9 &   &   &   &   &   &   &   & 6 & 6\rlap{$_?$} \\
          10 &  &   &   &   &   &   &   &   & 6\rlap{$_?$} 
        \end{tabular}
        \hfil
      \begin{tabular}{c|cccccccccccccccccccccccccccc}
          & 2 & 3 & 4 & 5 & 6 & 7 & 8 & 9 & 10 & 11 & 12 & 13 & 14 \\\hline
        2 & 1 & 1 & 1 & 1 & 1 & 1 & 1 & 1 & 1 & 1 & 1 & 1 & 1 \\
        3 &   & 1 & 2 & 2 & 2 & 2 & 2 & 2 & 2 & 2 & 2 & 2 & 2 \\
        4 &   &   & 2 & 2 & 3 & 3 & 3 & 3 & 3 & 3 & 3 & 3 & 3 \\
        5 &   &   &   & 2 & 3 & 3 & 4 & 4 & 4 & 4 & 4 & 4 & 4 \\
        6 &   &   &   &   & 3 & 3 & 4 & 4 & 5 & 5 & 5 & 5 & 5 \\
        7 &   &   &   &   &   & 3 & 4 & 4 & 5 & 5 & 6 & 6 & 6 \\
        8 &   &   &   &   &   &   & 4 & 4 & 5 & 5 & 6 & 6 & 7 \\
        9 &   &   &   &   &   &   &   & 4 & 5 & 5 & 6 & 6 & 7 \\
        10 &   &   &   &   &   &   &   &   & 5 & 5 & 6 & 6 & 7 \\
        11 &   &   &   &   &   &   &   &   &   & 5 & 6 & 6 & 7 \\
        12 &   &   &   &   &   &   &   &   &   &   & 6 & 6 & 7 \\
        13 &   &   &   &   &   &   &   &   &   &   &   & 6 & 7 \\
        14 &   &   &   &   &   &   &   &   &   &   &   &   & 7 
      \end{tabular}
    \end{center}

   The table on the right lists the numbers for $A_n$ and~$S_n$, which also turn out to
   be the same. For these groups, the invariant constraints make the problem easier, so
   that we were able to cover slightly larger values of $n$ and~$u$. All given numbers
   have been proved to be optimal. Note that a regular pattern emerges: we seem to have
   the formula $k=\min(u-1,\lfloor n/2\rfloor)$.

\end{document}